\documentstyle[twoside]{paper}

\begin{document}
\begin{center}
{\Huge Biorthogonal polynomials and 

\vspace{18pt}
zero-mapping transformations}

\vspace{32pt}
{\Large Arieh Iserles\footnote{Department of Applied Mathematics and
Theoretical Physics, University of Cambridge, England.} and Syvert P.
N\o rsett\footnote{Institute of Mathematical Sciences, Norwegian
Institute of Technology, Trondheim, Norway.}}

\vspace{36pt}
\parbox[t]{130mm}{{\bf Abstract} The authors have presented in
\cite{IN2} a technique to generate transformations $\cal T$ of the set
${\Bbb P}_n$ of $n$th degree polynomials to itself such that if
$p\in{\Bbb P}_n$ has all its zeros in $(c,d)$ then ${\cal T}\{p\}$ has
all its zeros in $(a,b)$, where $(a,b)$ and $(c,d)$ are given real
intervals. The technique rests upon the derivation of an explicit form
of biorthogonal polynomials whose Borel measure is strictly sign
consistent and such that the ratio of consecutive generalized moments is a
rational $[1/1]$ function of the parameter. Specific instances of
strictly sign consistent measures that have been debated in \cite{IN2}
include $x^\mu\D\psi(x)$, $\mu^x\D\psi(x)$ and $x^{\log_q\mu}\D\psi(x)$,
$q\in(0,1)$. In this paper we identify all measures $\psi$ such that
their consecutive generalized moments have a rational $[1/1]$
quotient, thereby characterizing all possible zero-mapping
transformations of this kind.}

\end{center}

\vfill
\eject

\section{Zero-mapping transformations}

In the present paper we wish to return to a theme that has been
already deliberated in \cite{IN2,INS1}. Let $(a,b)$ and $(c,d)$ be two
nonempty real intervals and denote by ${\Bbb P}_n$ the set of $n$th
degree polynomials. We are interested in linear transformations ${\cal
T}: {\Bbb P}_n \rightarrow {\Bbb P}_n$ such that any polynomial with
all its zeros in $(c,d)$ is mapped into a polynomial with all its
zeros in $(a,b)$.

A trivial instance of such a transformation is
$${\cal T}\left\{\sum_{k=0}^n a_k x^k\right\} =\sum_{k=0}^n ka_k
x^k,$$ 
which maps real zeros into real zeros. Another example, with
ubiquitous applications, is
$${\cal T}\left\{\sum_{k=0}^n a_k x^k\right\} =\sum_{k=0}^n
\Frac{1}{k!} a_k x^k$$
-- it maps positive zeros to positive zeros. Both follow from the
theory of {\em multiplier sequences\/}, well known since the pioneering
work of E. Laguerre \cite{L1} and of G. P\'olya and I. Schur
\cite{PS1}. The transformation
$${\cal T}\left\{\sum_{k=0}^n a_k x^k\right\} =\sum_{k=0}^n a_k
T_k(x),$$
where $T_k$ stands for the $k$th Chebyshev polynomial,
maps positive zeros to positive zeros -- the proof is elementary
(although perhaps surprising), takes just few lines of undergraduate
mathematics and we challenge the reader to find it without reference
to \cite{IS1}.

More examples of `zero-mapping' transformations are available in
literature, e.g.\ \cite{ASI1,P1}. A powerful technique for the
generation of such constructs has been presented by the authors in
\cite{IN2}. For every $\mu\in(c,d)$ we let $\D\varphi(x,\mu)$ be a
Borel measure supported by $x\in(a,b)$. Following \cite{IN1}, we say
that $p\in{\Bbb P}_n \times C\left[(c,d)^n\right]$, $p\not\equiv0$, is a {\em
biorthogonal polynomial\/} if
$$\int_a^b p(x;\mu_1,\ldots,\mu_n)\D\varphi(x,\mu_\ell)=0,\qquad
\ell=1,2,\ldots,n.$$
Biorthogonal polynomials exist and are unique (up to a nonzero
multiplicative constant) if and only if $\D\varphi$ is {\em regular\/}
\cite{IN1}, that is, for distinct $\mu_1,\mu_2,\ldots,\mu_n$,
$$\det\left[\begin{array}{llll}
\int_a^b \D\varphi(x,\mu_1) & \int_a^b \D\varphi(x,\mu_2) & \cdots &
\int_a^b \D\varphi(x,\mu_n)\\
\int_a^b x\D\varphi(x,\mu_1) & \int_a^b x\D\varphi(x,\mu_2) & \cdots &
\int_a^b x\D\varphi(x,\mu_n)\\
\C{\vdots} & \C{\vdots} & & \C{\vdots}\\
\int_a^b x^{n-1}\D\varphi(x,\mu_1) & \int_a^b
x^{n-1}\D\varphi(x,\mu_2) & \cdots & \int_a^b x^{n-1} \D\varphi(x,\mu_n)
	    \end{array}\right]\neq0,\quad n=1,2,\ldots.$$
We henceforth assume regularity.

Provided that $\varphi$ is a {\em strictly sign consistent (SSC)\/}
function, i.e.\ that
$$\det\left[\begin{array}{cccc}
\varphi(x_1,\mu_1) & \varphi(x_1,\mu_2) & \cdots &
\varphi(x_1,\mu_n)\\
\varphi(x_2,\mu_1) & \varphi(x_2,\mu_2) & \cdots &
\varphi(x_2,\mu_n)\\
\vdots & \vdots & & \vdots\\
\varphi(x_n,\mu_1) & \varphi(x_n,\mu_2) & \cdots &
\varphi(x_n,\mu_n)\end{array} \right]\neq0$$
for every $a<x_1<x_2<\cdots<x_n<b$, $c<\mu_1<\mu_2<\cdots<\mu_n<d$, it
can be proved that all the zeros of $p_n$ reside in the set $(a,b)$
\cite{IN1}. Likewise, zeros of $p_n$ live in $(a,b)$ if
$\D\varphi(x,\mu)=\omega(x,\mu)\D\psi(x)$, where $\D\psi(x)$ is a
Borel measure and the function $\omega$ is SSC. Biorthogonal
polynomials -- and, with greater generality, biorthogonal functions --
feature in a wide variety of interesting applications, mainly in
numerical analysis \cite{B1,IN3}.

Suppose that $\varphi$ (or $\omega$) is indeed an SSC function. Given any
$\mu_1,\ldots,\mu_n\in(c,d)$, we define
\begin{equation}\label{0.1}
{\cal T}\left\{ \prod_{\ell=1}^n (x-\mu_\ell)\right\}
=p_n(x;\mu_1,\mu_2,\ldots,\mu_n).
\end{equation}
Given that the range of the linear operator $\cal T$ in \R{0.1} can be
extended to all of 
${\Bbb P}_n$, it is a zero-mapping transformation. Specifically, it
maps a polynomial with all its zeros in $(c,d)$ into a polynomial with
all its zeros in $(a,b)$.

For the transformation \R{0.1} to be of any interest, beyond the most
formal, we need to know {\em explicitly\/} the form of
$p_n(x;\mu_1,\mu_2,\ldots,\mu_n)$. Fortunately, it is demonstrated in
\cite{IN2} that $p_n$ can be described in a closed form in an
important special case. Thus, let $\vec{\rho}:=\{\rho_k\}_{k\in\Zp}$
be an infinite sequence of monic polynomials, such that
\begin{eqnarray*}
\rho_0(x)&\equiv&1,\\
\rho_1(x)&=&(x-\sigma_1),\\
\rho_2(x)&=&(x-\sigma_1)(x-\sigma_2)\\
&\vdots&\\
\rho_k(x)&=&(x-\sigma_1)(x-\sigma_2)\cdots(x-\sigma_k)=
\rho_{k-1}(x)(x-\sigma_k),\qquad k=1,2,\ldots, 
\end{eqnarray*}
and consider the {\em generalized moments\/} 
$$I_k(\mu)=\int_a^b \rho_k(x)\D\varphi(x,\mu),\qquad k=0,1,\ldots.$$
Note that $\D\varphi$ is regular if and only if
\begin{equation}\label{0.2}
\det\left[\begin{array}{llll}
I_0(\mu_1) & I_0(\mu_2) & \cdots & I_0(\mu_n)\\
I_1(\mu_1) & I_1(\mu_2) & \cdots & I_1(\mu_n)\\
\C{\vdots} & \C{\vdots} & & \C{\vdots}\\
I_{n-1}(\mu_1) & I_{n-1}(\mu_2) & \cdots & I_{n-1}(\mu_n)
	    \end{array}\right]\neq0,\qquad n=1,2,\ldots,
\end{equation}
for all $c<\mu_1<\mu_2<\cdots<\mu_n<d$.

If the sequence $\{I_{k+1}(\mu)/I_k(\mu)\}_{k\in\Zp}$ consists of
rational $[1/1]$ functions, i.e.\ 
\begin{equation}\label{0.3}
\frac{I_{k+1}(\mu)}{I_k(\mu)}=\frac{g_k(\mu)}{h_k(\mu)},
\end{equation}
where $g_k,h_k\in{\Bbb P}_1$ then, subject to $g_kh_k'-g_k'h_k\neq0$,
$x\in(c,d)$, $k\in\Zp$,
$$p_n(x;\mu_1,\mu_2,\ldots,\mu_n)=\sum_{k=0}^n d_k \rho_k(x),$$
where
$$\prod_{\ell=1}^n (x-\mu_\ell)=\sum_{k=0}^n d_k \prod_{j=0}^{k-1}
g_j(x) \prod_{j=k}^{n-1} h_j(x).$$
In other words, \R{0.1} becomes \begin{equation}\label{0.4}
T\left\{\sum_{k=0}^m d_k \prod_{j=0}^{k-1} g_j(x) \prod_{j=k}^{m-1}
h_j(x) \right\}=\sum_{k=0}^m d_k \rho_k(x).
\end{equation}
In particular, if $\omega$, say, is SSC then \R{0.4} maps polynomials
with all their zeros in $(c,d)$ into polynomials with all their zeros
in $(a,b)$ \cite{IN2}.

Fifteen examples of transformation of this form have been presented in
\cite{IN2}. For example, letting $\D\varphi(x,\mu)=(\Gamma(\mu))^{-1}
x^{\mu-1}\E^{-x}\D x$, $(a,b)=(c,d)=(0,\infty)$, and $\sigma_k\equiv0$, we
have $I_k(\mu)=(\mu)_k$, $k\in\Zp$. Here $(z)_k$ is the {\em
Pochhammer symbol\/} \cite{R1},
$$(z)_0=1\qquad\mbox{and}\qquad
(z)_k=(z)_{k-1}(z+k-1)=\prod_{\ell=0}^{k-1} (z+\ell),\quad k=1,2,\ldots.$$
Therefore $g_k(\mu)=k+\mu$, $h_k(\mu)\equiv1$,
$g_kh_k'-g_kh_k'\equiv1$ and, $x^\mu$ being SSC \cite{IN1}, it
follows from \R{0.4} that the {\em Laguerre transformation\/}
$${\cal T}\left\{ \sum_{k=0}^n d_k (x)_k\right\} =\sum_{k=0}^n d_k
x^k$$
maps polynomials with positive zeros into polynomials with positive
zeros.

All the transformations in \cite{IN2} follow from eight strictly sign
consistent choices of
$\omega$, namely (a)~$x^\mu$, $x,\mu>0$, (b)~$\mu^x$, $x\in\RR$, $\mu>0$,
(c)~$x^{\log_q\mu}$, $x,\mu>0$, $q\in(0,1)$, (d)~$x^{\log_q\mu}$,
$x,\mu>0$, $q>1$, (e)~$\Gamma(x+\mu)$, $x,\mu>0$,
(f)~$1/\Gamma(x+\mu)$, $x,\mu>0$, (g)~$1/(x\mu;q)_\infty$,
$x,\mu,q\in(0,1)$, (h)~$(-x\mu;q)_\infty$, $x>0$, $\mu,q\in(0,1)$. The
notation $(z;q)_k$ stands for the {\em Gau\ss--Heine symbol\/}
\cite{S1},
$$(z;q)_0=1\qquad\mbox{and}\qquad
(z;q)_k=(z;q)_{k-1}(1-q^{k-1}z)=\prod_{\ell=0}^{k-1}(1-q^\ell z),\quad
k=1,2,\ldots,\infty.$$
In other words, each transformation is obtained by choosing one of the
above functions $\omega$, in tandem with a specific choice of $\D\psi$
and $\vec{\rho}$ which is consistent with \R{0.3}.

In the present paper we adopt a complementary approach. Thus, given
$\omega$, we attempt to identify {\em all\/} Borel measures $\D\psi$
and sets $\vec{\rho}$ such that \R{0.3} is true. Specifically, we
consider the three choices (a)--(c). There are in \cite{IN2} seven
transformations corresponding to these choices. We prove in the sequel
that, up to linear mapping of $x$ and $\mu$, this almost
exhausts the list of all possible transformations -- just a single
transformation has been missed in \cite{IN2}!

We expect to return to this issue in a future paper, characterizing
all transformations associated with the choices (d)--(h) and
consistent with \R{0.3}.

\section{The measure $\D\varphi(x,\mu)=x^\mu\D\psi(x)$, $\mu>0$}

Let $E=(a,b)\subseteq(0,\infty)$ be the least real interval that
contains the essential support of $\D\psi$,
$\{\sigma_\ell\}_{\ell=1}^\infty$ be a set of arbitrary real numbers
and define $\rho_k(x)=\prod_{\ell=1}^k (x-\sigma_\ell)$,
$k=0,1,\ldots$. The {\em moments\/} of
$\vec{\rho}=\{\rho_k\}_{k=0}^\infty$ are
$$I_k(\mu,\vec{\rho}):=\int_E x^\mu \rho_k(x)\D\psi(x),\qquad
k=0,1,\ldots,$$
therefore 
\begin{eqnarray}
I_k(\mu,\vec{\rho})&=&\int_E
x^{\mu-1}(x-\sigma_{k+1}+\sigma_{k+1})\rho_k(x)\D\psi(x) \label{1.1}\\
&=&I_{k+1}(\mu-1,\vec{\rho}) +\sigma_{k+1}I_k(\mu-1,\vec{\rho}),\qquad
k=0,1,\ldots. \nonumber
\end{eqnarray}

The following trivial observation will be repeatedly used in the
sequel.

\vspace{8pt}
\noindent {\bf Proposition 1} Suppose that there exists $k\in\Zp$ such
that
\begin{equation}\label{1.1A}
\frac{I_{k+1}(\mu)}{I_k(\mu)}\equiv{\rm const.}
\end{equation}
Then $\D\varphi$ is not regular.

{\em Proof.} Follows at once from \R{0.2}, since \R{1.1A} implies that
two columns in the matrix are proportional and the determinant must
therefore vanish for $k\geq n$. \QED

\vspace{8pt}
We wish to characterize all $\D\psi$ and $\{\sigma_\ell\}$ so that
\begin{equation}\label{1.2}
\frac{I_{k+1}(\mu,\vec{\rho})}{I_k(\mu,\vec{\rho})}
=\frac{\alpha_k+\beta_k \mu}{\gamma_k+\delta_k\mu},\qquad k=0,1,\ldots,
\end{equation}
for arbitrary real constants $\alpha_k,\beta_k,\gamma_k,\delta_k$ such
that $|\alpha_k|+|\beta_k|, |\gamma_k|+|\delta_k|>0$, $k=0,1,\ldots$.

Substituting \R{1.1} into \R{1.2}, we have for all $k=0,1,\ldots$
\begin{eqnarray*}
\frac{\alpha_k+\beta_k \mu}{\gamma_k+\delta_k\mu}&=&
\frac{I_{k+2}(\mu-1,\vec{\rho})+\sigma_{k+2}
I_{k+1}(\mu-1,\vec{\rho})}{I_{k+1}(\mu-1,\vec{\rho}) +\sigma_{k+1}
I_k(\mu-1,\vec{\rho})}\\
&=&\frac{\displaystyle
\frac{I_{k+2}(\mu-1,\vec{\rho})}{I_{k+1}(\mu-1,\vec{\rho})}+ 
\sigma_{k+2}}{\displaystyle
1+\sigma_{k+1}\frac{I_k(\mu-1,\vec{\rho})}{I_{k+1}(\mu-1,
\vec{\rho})}} =\frac{\displaystyle
\frac{\alpha_{k+1}+\beta_{k+1}(\mu-1)}
{\gamma_{k+1}+\delta_{k+1}(\mu-1)} +\sigma_{k+2}} {\displaystyle 1
+\sigma_{k+1}
\frac{\gamma_k+\delta_k(\mu-1)}{\alpha_k+\beta_k(\mu-1)}} \\
&=&\frac{(\alpha_k-\beta_k)+\beta_k\mu}{(\gamma_{k+1}-\delta_{k+1})
+\delta_{k+1}\mu}\\
&&\quad\mbox{}\times
\frac{(\alpha_{k+1}-\beta_{k+1}+\sigma_{k+2}(\gamma_{k+1}-\delta_{k+1}))
+(\beta_{k+1}+\sigma_{k+2}\delta_{k+1})\mu}
{(\alpha_k-\beta_k+\sigma_{k+1}(\gamma_k-\delta_k))
+(\beta_k+\sigma_{k+1}\delta_k)\mu}.
\end{eqnarray*}
We shift $\mu\rightarrow(\mu+1)$, thus obtaining the cubic identity
\begin{eqnarray}
&&[(\alpha_k+\beta_k)+\beta_k\mu][\gamma_{k+1}+\delta_{k+1}\mu]
[(\alpha_k+\sigma_{k+1}\gamma_k) +(\beta_k+\sigma_{k+1} \delta_k) \mu]
\label{1.3}\\
&=&[\alpha_k+\beta_k\mu][(\gamma_k+\delta_k)+\delta_k\mu]
[(\alpha_{k+1}+\sigma_{k+2}\gamma_{k+1})
+(\beta_{k+1}+\sigma_{k+2}\delta_{k+1})\mu]. \nonumber
\end{eqnarray}
We distinguish between the following cases:

\vspace{8pt}
\noindent {\bf Case I:} $(\alpha_k+\beta_k)+\beta_k\mu$ is a constant
multiple of $\alpha_k+\beta_k \mu$.

Hence $\beta_k=0$ and, without loss of generality, $\alpha_k=1$.
Therefore \R{1.3} reduces to
\begin{eqnarray}
&&[\gamma_{k+1}+\delta_{k+1}\mu][(1+\sigma_{k+1}\gamma_k)+\sigma_{k+1}\delta_k
\mu] \label{1.3a}\\ 
&=&[(\gamma_k+\delta_k)+\delta_k\mu][(\alpha_{k+1}+\sigma_{k+2}\gamma_{k+1})
+(\beta_{k+1}+\sigma_{k+2}\delta_{k+1})\mu]. \nonumber
\end{eqnarray}
There are thus two possible subcases:

\vspace{6pt}
\noindent {\bf Subcase I.1:} $\gamma_{k+1}+\delta_{k+1}\mu=C(\gamma_k+
\delta_k+\delta_k\mu)$ for some $C\neq0$.

We thus deduce that
\begin{eqnarray*}
\alpha_{k+1}&=&C(1+\sigma_{k+1}\gamma_k-\sigma_{k+2}\gamma_{k+1}),\\
\beta_{k+1}&=&C(\sigma_{k+1}\delta_k-\sigma_{k+2}\delta_{k+1}),\\
\gamma_{k+1}&=&C(\gamma_k+\delta_k),\\
\delta_{k+1}&=&C\delta_k.
\end{eqnarray*}
However, since the numerator and the denominator of a rational
function can be rescaled by a nonzero constant, we may assume without
loss of generality that $C=1$. Therefore
\begin{equation}\label{1.4}
\left. \begin{array}{rcl}
\beta_k &=& 0,\\
\alpha_{k+1} &=&
1+(\sigma_{k+1}-\sigma_{k+2})\gamma_k-\sigma_{k+2}\delta_k,\\ 
\beta_{k+1} &=& (\sigma_{k+1}-\sigma_{k+2})\delta_k,\\
\gamma_{k+1} &=& \gamma_k+\delta_k,\\
\delta_{k+1} &=& \delta_k. \end{array}\qquad\qquad\right\}
\end{equation}

\vspace{6pt}
\noindent{\bf Subcase I.2:} $1+\sigma_{k+1}\gamma_k+\sigma_{k+1}
\delta_k \mu=C(\gamma_k+\delta_k +\delta_k\mu)$, $C\neq0$.

Thus, either $\delta_k=0$ or $C=\sigma_{k+1}\neq0$. In the first case
$I_{k+1}/I_k$ is constant as a function of $\mu$ and regularity is
lost. In the second case \R{1.3a} yields
\begin{eqnarray*}
\alpha_{k+1}&=&(\sigma_{k+1}-\sigma_{k+2})\gamma_{k+1},\\
\beta_{k+1}&=&(\sigma_{k+1}-\sigma_{k+2})\delta_{k+1}.
\end{eqnarray*}
Since $|\alpha_{k+1}|+|\beta_{k+1}|\neq0$, we deduce that
$\sigma_{k+2}\neq \sigma_{k+1}$ and $I_{k+2}/I_{k+1}$ is a constant.
This, again, contradicts regularity.

Thus, we deduce that {\bf Case I} necessarily implies \R{1.4}

\vspace{8pt}
\noindent{\bf Case II:} $(\alpha_k+\beta_k)+\beta_k\mu$ is a constant
nonzero multiple of $(\alpha_{k+1}+\sigma_{k+2}\gamma_{k+1})+
(\beta_{k+1}+\sigma_{k+2}\delta_{k+1})\mu$.

Thus, there exists $C_1\neq0$ such that
\begin{eqnarray*}
\alpha_{k+1}&=&C_1(\alpha_k+\beta_k)-\sigma_{k+2}\gamma_{k+1},\\
\beta_{k+1}&=&C_1\beta_k-\sigma_{k+2}\delta_{k+1}.
\end{eqnarray*}
Again, there are two possibilities in \R{1.3}. 

\vspace{6pt}
\noindent {\bf Subcase II.1:} There exists $C_2\neq0$ such that
\begin{eqnarray*}
\gamma_{k+1}+\delta_{k+1}\mu&=&C_2(\gamma_k+\delta_k+\delta_k\mu),\\
C_2((\alpha_k+\sigma_{k+1}\gamma_k)+(\beta_k+\sigma_{k+1}\delta_k)
\mu) &=&C_1(\alpha_k+\beta_k\mu).
\end{eqnarray*}

Therefore $\delta_{k+1}=C_2\delta_k$,
$\gamma_{k+1}=C_2(\gamma_k+\delta_k)$ and
\begin{eqnarray*}
(C_2-C_1)\alpha_k&=&-C_2\sigma_{k+1}\gamma_k,\\
(C_2-C_1)\beta_k&=&-C_2\sigma_{k+1}\delta_k.
\end{eqnarray*}

If $C_1\neq C_2$ then, again, $I_{k+1}/I_k$ is a constant and
regularity is lost. Hence necessarily $C_2=C_1$ and, since
$\gamma_k=\delta_k=0$ is impossible, we have $\sigma_{k+1}=0$. Thus,
we may assume without loss of generality that $C_1=C_2=1$ and obtain
\begin{equation}\label{1.5}
\left. \begin{array}{rcl}
\sigma_{k+1}&=&0,\\
\alpha_{k+1}&=&\alpha_k+\beta_k-\sigma_{k+2}(\gamma_k+\delta_k),\\
\beta_{k+1}&=&\beta_k-\sigma_{k+2}\delta_k,\\
\gamma_{k+1}&=&\gamma_k+\delta_k,\\
\delta_{k+1}&=&\delta_k.
       \end{array}\right\}
\end{equation}

\vspace{6pt}
\noindent{\bf Subcase II.2:} $C_2\neq0$ exists so that
\begin{eqnarray*}
\gamma_{k+1}+\delta_{k+1}\mu&=&C_2(\alpha_k+\beta_k\mu),\\
C_2((\alpha_k+\sigma_{k+1}\gamma_k)+(\beta_k+\sigma_{k+1}\delta_k) \mu)
&=&C_1(\gamma_k +\delta_k+\delta_k\mu).
\end{eqnarray*}

We deduce that $\gamma_{k+1}=C_2\alpha_k$, $\delta_{k+1}=C_2\beta_k$
and
\begin{eqnarray*}
\alpha_{k+1}&=&(C_1-\sigma_{k+2}C_2)\alpha_k,\\
\beta_{k+1}&=&(C_1-\sigma_{k+2}C_2)\beta_k.
\end{eqnarray*}
Therefore $I_{k+2}/I_{k+1}$ is a constant multiple of $I_{k+1}/I_k$,
and this cannot coexist with regularity -- the proof is identical to
that of Proposition 1. Hence, this subcase is impossible.

\vspace{8pt}
\noindent {\bf Case III:} $(\alpha_k+\beta_k)+\beta_k\mu$ is a
constant multiple of $(\gamma_k+\delta_k)+\delta_k\mu$.

It is obvious in that case that $I_{k+1}/I_k$ is a constant and this
is ruled out by regularity.

\vspace{8pt}
The above three cases exhaust all possibilities. Therefore, we deduce
that for every $k=0,1,\ldots$ either \R{1.4} or \R{1.5} must hold.
This, in particular, implies that
\begin{equation}\label{1.6}
\delta_k\equiv\delta_0,\quad\gamma_k=\gamma_0+k\delta_0,\qquad k=0,1,\ldots.
\end{equation}

Suppose first that $\delta_0=0$. Hence, without loss of generality,
$\gamma_k\equiv1$ and 
$$\frac{I_{k+1}(\mu)}{I_k(\mu)}=\alpha_k+\beta_k\mu.$$
Regularity thus requires $\beta_k\neq0$ and this rules out \R{1.4}. We
deduce that $\delta_0=0$ implies \R{1.5} for all $k=0,1,\ldots$. This
results in the explicit form
\begin{equation}\label{1.7}
\left.\begin{array}{rcl}
\sigma_{k+1}&=&0,\\
\alpha_k&=&\alpha_0+k\beta_0,\\
\beta_k&\equiv&\beta_0,\\
\gamma_k&\equiv&1,\\
\delta_k&\equiv&0.
\end{array}\right\}
\end{equation}

Next we consider the case $\delta_0\neq0$ and assume without loss of
generality that $\delta_0=1$. Either \R{1.4} or \R{1.5} must
hold for each $k\in\Zp$ and we commence by assuming that
integers
$$0=m_0<n_0<m_1<n_1<\cdots$$
exist so that for all $\ell=0,1,\ldots$
\begin{eqnarray*}
k\in\{m_\ell,m_\ell+1,\ldots,n_\ell-1\} &\qquad\Rightarrow\qquad&
\mbox{\R{1.4}},\\
k\in\{n_\ell,n_\ell+1,\ldots,m_{\ell+1}-1\} &\qquad\Rightarrow\qquad&
\mbox{\R{1.5}}.
\end{eqnarray*}
Since $\beta_{k+1}=0$ implies $\sigma_{k+1}=\sigma_{k+2}$ in \R{1.4},
we deduce that
\begin{equation}\label{1.8}
\sigma_k=\sigma_{m_\ell+1},\qquad k=m_\ell+1,\ldots,n_\ell.
\end{equation}
By applying similar argument to \R{1.5} we deduce
\begin{equation}\label{1.9}
\sigma_k=0,\qquad k=n_\ell,n_\ell+1,\ldots,m_{\ell+1}.
\end{equation}
Thus, \R{1.5} implies that $\beta_k=\beta_{n_\ell}$,
$k=n_\ell,n_\ell+1,\ldots,m_{\ell+1}-1$ and, to obtain
$\beta_{m_{\ell+1}}=0$, we need
$\beta_{n_\ell}=\sigma_{m_{\ell+1}+1}$. But
$\beta_{n_\ell}=\sigma_{m_\ell}$ and we deduce that
$\sigma_{m_\ell}\equiv \sigma$, say, for all $\ell=0,1,\ldots$.
Therefore, letting $g_k(\mu)=\alpha_k+\beta_k\mu$,
\begin{eqnarray*}
k=m_\ell,\ldots,n_\ell-1:&\qquad& g_k(\mu)\equiv\alpha_0-\ell\sigma,\\
k=n_\ell,\ldots,m_{\ell+1}-1: &\qquad&
g_k(\mu)=\alpha+\sigma(\gamma_0+k-\ell-1) +\sigma\mu.
\end{eqnarray*}
This formula is consistent with \R{1.2}, but we need to check whether
it is also consistent with the fact that $\{I_k(\mu)\}_{k=0}^\infty$
corresponds to a Hamburger moment sequence for every $\mu>0$. Let
$$\tilde{I}_k(\mu):=\int_E x^{k+\mu}
\D\psi(x)=I_k(\mu,\{x^j\}_{j=0}^\infty),\qquad k=0,1,\ldots.$$
Then, by \cite{Ak}, we need
$$\Delta_{n}:=\det\left[ \begin{array}{llll}
\tilde{I}_0(\mu) & \tilde{I}_{1}(\mu) & \cdots &
\tilde{I}_{n}(\mu)\\ 
\tilde{I}_{1}(\mu) & \tilde{I}_{2}(\mu) & \cdots &
\tilde{I}_{n+1}(\mu),\\ 
\C{\vdots} & \C{\vdots} & & \C{\vdots}\\
\tilde{I}_{n}(\mu) & \tilde{I}_{n+1}(\mu) & \cdots & \tilde{I}_{2n}(\mu)
	     \end{array}\right]>0,\qquad n=0,1,\ldots,\quad \mu>0.$$

Let us first assume that $n_0\geq 2$. Then, letting
$h_k(\mu)=\gamma_k+\delta_k\mu$, and bearing in mind that
$I_k=I_k(\,\cdot\,,\{(x-\sigma)^\ell\}_{\ell=0}^\infty)$, 
\begin{eqnarray*}
I_1&=&\frac{g_0}{h_0}I_0,\quad I_2=\frac{g_0g_1}{h_0h_1}I_0,\\
\tilde{I}_0&=&I_0,\quad \tilde{I}_1=I_1+\sigma I_0,\quad
\tilde{I}_2=I_2+2\sigma I_1+\sigma^2 I_0
\end{eqnarray*}
imply 
$$\tilde{I}_1=\left(\frac{\alpha_0}{\gamma_0+\mu}+\sigma\right)\tilde{I}_0$$
and 
$$\tilde{I}_2=\left(\frac{\alpha_0^2}{(\gamma_0+\mu)(\gamma_0
+\mu+1)}+2\frac{\sigma\alpha_0}{\gamma_0+\mu}+\sigma^2\right)\tilde{I}_0,$$
therefore
$$\Delta_1=-\frac{\alpha_0^2}{(\gamma_0+\mu)^2(\gamma_0+\mu+1)}\tilde{I}_0^2
<0$$
for sufficiently large $\mu$. We deduce that $n_0\geq2$ is impossible.

The remaining case is $n_0=1$. $\tilde{I}_1$ remains intact, whereas
$\tilde{I}_2=I_2+\sigma I_1+\sigma^2 I_0$ and
$$\Delta_1=-\frac{\alpha_0(\alpha_0+\sigma(\gamma_0+\mu))}{(\gamma_0+\mu)^2
(\gamma_0+\mu+1)}\tilde{I}_0^2<0$$
for $\mu\gg0$. We deduce that $0=m_0<n_0<m_1<\cdots$ is impossible.

Finally, we check the case whereby there exist
$$0=n_0<m_1<n_1<m_2<\cdots$$
so that
\begin{eqnarray*}
k\in\{n_\ell,n_\ell+1,\ldots,m_{\ell+1}-1\} &\qquad\Rightarrow\qquad&
\mbox{\R{1.5}},\\
k\in\{m_\ell,m_\ell+1,\ldots,n_\ell-1\} &\qquad\Rightarrow\qquad&
\mbox{\R{1.4}}.
\end{eqnarray*}

\R{1.8} and \R{1.9} follow as before and, in addition,
\begin{eqnarray*}
k=n_\ell,\ldots,m_{\ell+1}-1: &\qquad&
g_k(\mu)=\alpha_0+\sigma(k-\ell)+\sigma\mu,\\
k=m_\ell,\ldots,n_\ell-1: &\qquad& g_k(\mu)=\alpha_0-\sigma(\gamma_0+\ell-1).
\end{eqnarray*}

Suppose that such $m_1\geq1$ exists and let
$$\D\varphi^*(x,\mu):=x^{m_1}\D\varphi(x,\mu).$$
Therefore
$$\tilde{I}^*_\ell(\mu):=\int_E x^\ell
\D\varphi^*(x,\mu)=\tilde{I}_{\ell+m_1}(\mu),\qquad \ell=0,1,\ldots,$$
and, proceeding as before (replacing $\tilde{I}_0$ by
$\tilde{I}_{m_1}$ etc.) we can prove that
$\{\tilde{I}^*_\ell\}_{\ell=0}^\infty$ cannot be a moment sequence.

We conclude that no such $m_1$ exists, hence, necessarily, \R{1.5}
{\em is true for all $k=0,1,\ldots$.\/} We obtain
\begin{equation}\label{1.10}
\left.\begin{array}{rcl}
\sigma_k&\equiv&0,\\
\alpha_k&=&\alpha_0+k\beta_0,\\
\beta_k&\equiv&\beta_0,\\
\gamma_k&=&\gamma_0+k\delta_0,\\
\delta_k&\equiv&\delta_0.
      \end{array}\right\}
\end{equation}
Note that \R{1.7} is a special case.

If $\delta_0=0$ we obtain (with $\gamma_0=1$)
$$g_k(\mu)=\alpha_0+k\beta_0+\beta_0\mu,\quad h_k(\mu)\equiv1,\qquad
k=0,1,\ldots$$
and this, in the terminology of \cite{IN2}, is the {\em Laguerre
transformation.\/} On the other hand, if $\delta_0\neq0$ then, letting
$\delta_0=1$,
$$g_k(\mu)=\alpha_0+k\beta_0+\beta_0\mu,\quad
h_k(\mu)=\gamma_0+k+\mu,\qquad k=0,1,\ldots,$$
namely the {\em Jacobi transformation\/} \cite{IN2}.

\vspace{8pt}
\noindent {\bf Theorem 2} The only transformations consistent with the
asserted form of $\D\varphi$ are the Laguerre transformation and the
Jacobi transformation. \QED

\section{The measure $\D\varphi(x,\mu)=\mu^x\D\psi(x)$, $\mu>0$}

$E$ is now a measurable subset of $(-\infty,\infty)$, otherwise we
employ the notation from the previous section. Thus, differentiating
$$I_k(\mu,\vec{\rho})=\int_E \rho_k(x)\mu^x\D\psi(x),\qquad
k=0,1,\ldots,$$
with respect to $\mu$ we obtain
$$I_k'(\mu)=\int_E \rho_k(x)x\mu^{x-1}\D\psi(x)=\mu^{-1}\int_E
\rho_k(x) ((x-\sigma_{k+1})+\sigma_{k+1})\mu^x \D\psi(x)$$
(here and elsewhere we abbreviate $I_k(\mu)=I_k(\mu,\vec{\rho})$)
and deduce the equation
\begin{equation}\label{2.1}
\mu I_k'(\mu)=I_{k+1}(\mu)+\sigma_{k+1}I_k(\mu),\qquad k=0,1,\ldots.
\end{equation}
Assuming \R{1.2}, we hence obtain from \R{2.1} the linear differential
equation
\begin{equation}\label{2.2}
I_k'(\mu)=\frac{1}{\mu} \left\{
\sigma_{k+1}+\frac{\alpha_k+\beta_k\mu} {\gamma_k+\delta_k\mu}
\right\} I_k(\mu)
\end{equation}

We distinguish among the following cases:

\vspace{8pt}
\noindent{\bf Case I:} $\gamma_k=0$ and, without loss of generality,
$\delta_k=1$.

Therefore, by \R{2.2},
$$I_k'(\mu)=\left\{(\sigma_{k+1}+\beta_k)\frac{1}{\mu}+\frac{\alpha_k}{\mu^2}
\right\}I_k(\mu)$$
and integration yields the explicit form
$$I_k(\mu)=I_k(\mu_0)
\left(\frac{\mu}{\mu_0}\right)^{\sigma_{k+1}+\beta_k}
\exp\left(-\alpha_k(\mu^{-1} -\mu_0^{-1})\right),$$
where $\mu_0>0$. In particular, choosing without loss of generality
$\mu_0=1$, we obtain
\begin{equation}\label{2.3}
I_k(\mu)=I_k(1)\E^{\alpha_k} \mu^{\sigma_{k+1}+\beta_k} \E^{-\alpha_k/\mu}.
\end{equation}

\vspace{8pt}
\noindent{\bf Case II:} $\delta_k=0$ and, without loss of generality,
$\gamma_k=1$.

Again, we solve \R{2.2} explicitly with an initial condition at
$\mu_0=1$, say, to obtain
\begin{equation}\label{2.4}
I_k(\mu)=I_k(1)\E^{-\beta_k} \mu^{\sigma_{k+1}+\alpha_k} \E^{\beta_k\mu}.
\end{equation}

\vspace{8pt}
\noindent {\bf Case III:} $\gamma_k,\delta_k\neq0$ and, without loss
of generality, $\delta_k=1$.

Therefore
$$I_k'(\mu)=\left\{\left(\sigma_{k+1}+\frac{\alpha_k}{\gamma_k}\right)
\frac{1}{\mu} -\left(\frac{\alpha_k}{\gamma_k}-\beta_k\right)
\frac{1}{\mu+\gamma_k} \right\}I_k(\mu),$$
with the solution
\begin{equation}\label{2.5}
I_k(\mu)=I_k(1)(1+\gamma_k)^{\alpha_k/\gamma_k-\beta_k}
\mu^{\sigma_{k+1}
+\alpha_k/\gamma_k}(\mu+\gamma_k)^{-\alpha_k/\gamma_k+\beta_k}. 
\end{equation}

\vspace{8pt}
Let us suppose that we have case I for $k$ and case II for $k+1$.
Therefore
$$\frac{I_{k+1}(\mu)}{I_k(\mu)}=\frac{I_{k+1}(1)}{I_k(1)}
\E^{\alpha_k+\beta_{k+1}}
\mu^{\sigma_{k+1}-\sigma_{k+2}+\beta_k-\alpha_{k+1}}
\E^{-(\alpha_k/\mu+\beta_{k+1}\mu)}.$$
But, since $\gamma_k=0$, $\delta_k=1$, we have
$$\frac{I_{k+1}(\mu)}{I_k(\mu)}=\frac{\alpha_k}{\mu}+\beta_k$$
and, necessarily, $\alpha_k/\mu+\beta_{k+1}\mu\equiv{\rm const}$. We
deduce that $\alpha_k=0$, but this, in tandem with $\gamma_k=0$,
contradicts regularity.

Similar contradiction is obtained if case I follows case II.

Now, if case I is followed by case III then, to get rid of the exponential
term, we require $\alpha_k=0$, and this again contradicts regularity.
Similarly, case II cannot be followed by case III. Finally, if case
III is followed by case I (case II) then again we need to eliminate an
exponential term, this requires $\alpha_{k+1}=0$ ($\beta_{k+1}=0$) and
is inconsistent with regularity.

We deduce that if any of cases I--III holds for one $k$ then it must
hold for all $k=0,1,\ldots$.

\vspace{8pt}
\noindent {\bf Case I} for all $k=0,1,\ldots$:

Therefore
$$\frac{I_{k+1}(\mu)}{I_k(\mu)}=\frac{I_{k+1}(1)}{I_k(1)}
\mu^{\sigma_{k+2}-\sigma_{k+1}+\beta_{k+1}-\beta_k}
\E^{(\alpha_k-\alpha_{k+1})(\mu-1)} =\frac{\alpha_k}{\mu}+\beta_k,\qquad
k=0,1,\ldots.$$
We deduce that
$$\alpha_{k+1}=\alpha_k=\cdots=\alpha_0.$$
Moreover, there are exactly two possibilities. Either $\beta_k=0$ and
$\sigma_{k+2}-\sigma_{k+1}+\beta_{k+1}-\beta_k=-1$ or $\alpha_0=0$. In
the second case $I_{k+1}/I_k$ is a constant and regularity is
violated, hence necessarily $\beta_k=0$ for all $k=0,1,\ldots$. We
hence deduce that $\sigma_{k+2}-\sigma_{k+1}=-1$, thus
\begin{equation}\label{2.6}
\left. \begin{array}{rcl}
\sigma_k&=&\sigma_1-k+1,\\
\alpha_k&\equiv&\alpha_0\neq0,\\
\beta_k&\equiv&0.
       \end{array}\right\}
\end{equation}
In other words,
$$I_k(\mu)=\left(\frac{\alpha_0}{\mu}\right)^k I_0(\mu),\qquad
k=0,1,\ldots,$$
and $\rho_k(x)=(x-\sigma_1)_k$. Under a substitution $x\rightarrow
-x+\sigma_1$, $\mu\rightarrow \alpha_0/\mu$ this yields {\em precisely\/}
the Charlier transformation from \cite{IN2}.

\vspace{8pt}
\noindent {\bf Case II} for all $k=0,1,\ldots$:

Since
$$\frac{I_{k+1}(\mu)}{I_k(\mu)}=\frac{I_{k+1}(1)}{I_k(1)}
\mu^{\sigma_{k+2}-\sigma_{k+1}+\alpha_{k+1}-\alpha_k}
\E^{(\beta_{k+1}-\beta_k)(\mu-1)}=\alpha_k+\beta_k\mu,$$
we deduce that $\beta_{k+1}=\beta_k=\cdots=\beta_0$. Again, there are
two alternatives -- either $\beta_0=0$ and this contradicts regularity
(since $I_{k+1}/I_k$ becomes a constant) or $\alpha_k=0$ and
$\alpha_{k+1}=\sigma_{k+1}-\sigma_{k+2}+1$. We deduce that
$\alpha_k\equiv0$ and $\sigma_k=\sigma_1+k-1$. Therefore
\begin{equation}\label{2.7}
\left. \begin{array}{rcl}
\sigma_k&=&\sigma_1+k-1,\\
\alpha_k&\equiv&0,\\
\beta_k&\equiv&\beta_0\neq0.
       \end{array}\right\}
\end{equation}
We conclude that
$$I_k(\mu)=(\beta_0\mu)^kI_0(\mu),\qquad k=0,1,\ldots,$$
and $\rho_k(x)=(-1)^k (-x+\sigma_1)_k$. Therefore $x\rightarrow
x-\sigma_1$, $\mu\rightarrow \mu/\beta_0$ yields, again, the Charlier
transformation.

\vspace{8pt}
\noindent {\bf Case III} for all $k=0,1,\ldots$. 

We now have
$$\frac{I_{k+1}(\mu)}{I_k(\mu)}=C\mu^{\sigma_{k+2}-\sigma_{k+1}+\alpha_{k+1}
/\gamma_{k+1}-\alpha_k/\gamma_k}
\frac{(\mu+\gamma_k)^{\alpha_{k+1}/\gamma_{k+1} -\beta_k}}
{(\mu+\gamma_{k+1})^{\alpha_{k+1}/\gamma_{k+1}-\beta_{k+1}}}
=\frac{\alpha_k+\beta_k \mu}{\mu+\gamma_k},$$
where $C\neq0$. Regularity implies that $\alpha_j/\gamma_j\neq
\beta_j$ for all $j\in\Zp$ (otherwise $I_{k+1}/I_k$ is constant).
Therefore, letting $\mu=-\gamma_k$ in
$$C\mu^{\sigma_{k+2}-\sigma_{k+1}+\alpha_{k+1}/\gamma_{k+1}-\alpha_k
/\gamma_k} (\mu+\gamma_k)^{\alpha_k/\gamma_k-\beta_k+1}
=(\alpha_k+\beta_k\mu)
(\mu+\gamma_{k+1})^{\alpha_{k+1}/\gamma_{k+1}-\beta_{k+1}}$$
demonstrates at once that $\gamma_{k+1}=\gamma_k\equiv\gamma_0\neq0$
and that either $\alpha_k=0$ or $\beta_k=0$. Thus, either
\begin{equation}\label{2.8}
\alpha_k=0,\qquad
\frac{\alpha_{k+1}}{\gamma_0}=1-\sigma_{k+2}+\sigma_{k+1},\qquad
\beta_{k+1}-\sigma_{k+2}+\sigma_{k+1}
\end{equation}
or
\begin{equation}\label{2.9}
\frac{\alpha_{k+1}}{\gamma_0}=\frac{\alpha_k}{\gamma_0}+\sigma_{k+1}
-\sigma_{k+2}, \qquad \beta_k=0,\qquad
\beta_{k+1}=\sigma_{k+1}-\sigma_{k+2}-1.
\end{equation}

Suppose that integers
$$0=m_0<n_0<m_1<n_1<\cdots$$
exist so that \R{2.8} holds for all $k=m_\ell,\ldots,n_\ell-1$ and
\R{2.8} holds for $k=n_\ell,\ldots,m_{\ell+1}-1$. The range of $\ell$
might be infinite or finite -- in the latter case we assume that the
largest value of $m_\ell$ or $n_\ell$, as the case might be, is
$\infty$. Thus, for example, if \R{2.8} is satisfied by all $k\in\Zp$
then this corresponds to $m_0=0$, $n_0=\infty$.

There is no loss of generality in assuming that \R{2.8} holds for
$k=0$, otherwise we map $x\rightarrow-x$. Hence, only the present case
need be considered. Moreover, we may assume without loss of generality
that $\sigma_1=0$, otherwise we replace $\sigma_k$ by
$\sigma_k-\sigma_1$ and shift $x\rightarrow x+\sigma_1$.

Let
$$M_\ell=\sum_{j=0}^\ell m_j,\qquad N_\ell=\sum_{j=0}^\ell n_j.$$
Long and tedious calculation affirms that
\begin{eqnarray}
\left.\begin{array}{rcl}
\alpha_k&\equiv&0,\quad \beta_k=\beta_0-k,\quad \gamma_k\equiv0,\quad
\delta_k\equiv1,\\
\sigma_{k+1}&=&N_{\ell-1}-M_\ell+k,\\
\rho_k(x)&=&(-1)^{N_{\ell-1}-M_\ell+k} (-x)_{N_{\ell-1}-M_\ell+k}
(x-\beta_0)_{M_\ell-N_{\ell-1}}
\end{array}\hspace*{-8pt}\right\}&& \hspace*{-16pt}
k=m_\ell,\ldots,n_\ell-1;\qquad\quad \label{2.10}\\
\left.\begin{array}{rcl}
\alpha_k&=&(k-\beta_0)\gamma_0,\quad \beta_k\equiv0,\quad
\gamma_k\equiv\gamma_0,\quad \delta_k\equiv1,\\
\sigma_{k+1}&=&\beta_0+N_\ell-M_\ell-k,\\
\rho_k(x)&=&(-1)^{N_\ell-M_\ell}(-x)_{N_\ell-M_\ell}
(x-\beta_0)_{M_\ell-N_\ell+k} \hspace*{32pt}
      \end{array}\right\}&& \hspace*{-16pt}
k=n_\ell,\ldots,m_{\ell+1}-1 \qquad \label{2.11}
\end{eqnarray}
and
\begin{eqnarray}
I_k(\mu)&=&(-1)^{N_{\ell-1}-M_\ell+k}
\frac{(-\beta_0)_k}{(\gamma_0+\mu)^k} \mu^{N_{\ell-1}-M_\ell+k}
\gamma_0^{M_\ell-N_{\ell-1}},\quad k=m_\ell,\ldots,n_\ell;
\label{2.12}\\
I_k(\mu)&=&(-1)^{N_\ell-M_\ell} \frac{(-\beta_0)_k}{(\gamma_0+\mu)^k}
\mu^{N_\ell-M_\ell} \gamma_0^{M_\ell-N_\ell+k},\hspace*{50pt}
k=n_\ell,\ldots,m_{\ell+1}. \label{2.13}
\end{eqnarray}
Note that both definitions match at $k=m_\ell$ and $k=n_\ell$ and that
they are indeed consistent -- as they should -- with 
$$\frac{I_{k+1}(\mu)}{I_k(\mu)}=\frac{\alpha_k+\beta_k\mu}{\gamma_0+\mu},
\qquad k=0,1,\ldots.$$

We will now single out a measure $\D\varphi$ of the desired form that
results in the required value of $\{I_k\}_{k\in\Zp}$. Thus, we let 
$$\D\varphi(x,\mu)=\mu^x\D\psi(x),\qquad x>0,$$
where $\D\psi$ is an atomic measure with jumps of
$$(-1)^j \left(1+\frac{\mu}{\gamma_0}\right)^{-\beta_0}
\frac{(-\beta_0)_j}{j!\gamma_0^j}$$
at $j=0,1,\ldots$. In other words, for every appropriate
function $f$ we have
$$\int_0^\infty f(x)\D\varphi(x,\mu)=\sum_{j=0}^\infty (-1)^j
\frac{(-\beta_0)_jf(j)}{j!}\left(\frac{\mu}{\gamma_0}\right)^j.$$
In particular, substituting the values of $\rho_k$ from \R{2.10} and
\R{2.11} respectively, we obtain
\begin{eqnarray}
I_k(\mu)&=&\left(1+\frac{\mu}{\gamma_0}\right)^{-\beta_0}
\sum_{j=N_{\ell-1}-M_\ell+k}^\infty (-1)^j 
\frac{(j-\beta_0)_{M_\ell-N_{\ell-1}}(-\beta_0)_j}{(j-N_{\ell-1}+M_\ell-k)!}
\left(\frac{\mu}{\gamma_0}\right)^j,\nonumber\\
&&\hspace*{200pt} k=m_\ell,\ldots,n_\ell, \label{2.14}\\
I_k(\mu)&=&\left(1+\frac{\mu}{\gamma_0}\right)^{-\beta_0}
\sum_{j=N_\ell-M_\ell}^\infty (-1)^j 
\frac{(j-\beta_0)_{M_\ell-N_\ell+k} (-\beta_0)_j}{(j-N_\ell+M_\ell)!}
\left(\frac{\mu}{\gamma_0}\right)^j,\nonumber\\
&& \hspace*{200pt} k=n_\ell,\ldots,m_{\ell+1}. \label{2.15} 
\end{eqnarray}
Our goal is to demonstrate that \R{2.12} is identical with \R{2.14}
and \R{2.13} coincides with \R{2.15}.

We now prove that \R{2.14} is the same as \R{2.12}. Letting
$s=M_\ell-N_{\ell-1}>0$ and commencing from \R{2.14}, we have for
every $k=m_\ell,\ldots,n_\ell$ 
\begin{eqnarray*}
I_k(\mu)&=&\left(1+\frac{\mu}{\gamma_0}\right)^{-\beta_0}
\sum_{j=k-s}^\infty (-1)^j \frac{(-\beta_0)_{j+s}}{(j+s-k)!}
\left(\frac{\mu}{\gamma_0}\right)^j\\
&=&(-1)^{k-s}(-\beta_0)_k \left(\frac{\mu}{\gamma_0}\right)^{k-s}
\left(1+\frac{\mu}{\gamma_0}\right)^{-\beta_0} 
\sum_{j=0}^\infty (-1)^j \frac{(-\beta_0+k)_j}{j!}
\left(\frac{\mu}{\gamma_0}\right)^j \\
&=&(-1)^{k-s}(-\beta_0)_k \left(\frac{\mu}{\gamma_0}\right)^{k-s}
\left(1+\frac{\mu}{\gamma_0}\right)^{-\beta_0} \times
\left(1+\frac{\mu}{\gamma_0}\right)^{\beta_0-k}\\
&=&(-1)^{k-s}\frac{(-\beta_0)_k}{(\gamma_0+\mu)^k} \mu^{k-s}\gamma_0^s.
\end{eqnarray*}
Therefore, we recover \R{2.12}.

The proof of the coincidence of \R{2.13} and \R{2.15} is identical. We
thus deduce that $\D\varphi$ is of the stipulated form.\footnote{We
have not proved that this $\D\varphi$ is the {\em unique\/} measure
with this property. This is quite straightforward since it is easy to
prove that it is determinate (e.g.\ with the Carleman criterion \cite{C1}) but
irrelevant since transformations are determined by moment sequences.}

In order to identify $\D\varphi$ we note that, to be a Borel measure,
the jumps must be nonnegative for all $j\in\Zp$. There are
exactly two possibilities. Firstly, $\beta_0,\gamma_0<0$ imply that
$$(-1)^j\frac{(-\beta_0)_j}{\gamma_0^j}=\frac{(|\beta_0|)_j}
{|\gamma_0|^j}>0,\qquad j=0,1,\ldots.$$
In that case $\mu\in(0,|\gamma_0|)$ and we recover the {\em Meixner\/}
transformation \cite{IN2}. Secondly, $\beta_0=N$ is a nonnegative
integer and $\gamma_0$. Thus,
$$(-1)^j\frac{(-\beta_0)_j}{\gamma_0^j}=\left\{ \begin{array}{lcl}
\frac{N!}{(N-j)!\gamma_0^j}>0 & \qquad & :j=0,1,\ldots,N,\\
0 &\qquad & :j=N+1,N+2,\ldots. 	\end{array}\right.$$
This, in the terminology of \cite{IN2}, is the {\em Krawtchouk\/}
transformation.

\vspace{8pt}
\noindent{\bf Theorem 3}  The only transformations consistent with the
asserted form of $\D\varphi$ are the Charlier, Meixner and Krawtchouk
transformation. \QED

\section{The measure $\D\varphi(x,\mu)=x^{\log_q\mu}\D\psi(x)$,
$q\in(0,1)$}

At first glance, this case is equivalent to $\D\varphi(x,\nu)=x^\nu
\D\psi(x)$, subject to the transformation $\eta=\log_q\mu$. However,
in that case the quotient of $I_{k+1}$ and $I_k$ will cease to be a
rational $[1/1]$ function of $\mu$, thereby violating our construction.

We consider $\D\varphi(x,\mu)=x^{\log_q\mu}\D\psi(x)$, where
$q\in(0,1)$, subject to the condition that $I_k(\mu,\vec{\rho})$ is
well-defined and bounded at $\mu=0$. The last condition sounds strange
-- after all, $\log_q\mu$ becomes unbounded as $\mu=0$ -- but it is
not! Recall that, by our assumption in this paper,
$$I_{k+1}(\mu,\vec{\rho})=\frac{\alpha_k+\beta_k\mu}{\gamma_k+
\delta_k\mu} I_k(\mu,\vec{\rho}),\qquad k=0,1,\ldots,$$
hence we require that $\gamma_k\neq0$, $k\in\Zp$, and that
$I_0(0,\vec{\rho})$ is bounded. However,
$$I_0(\mu,\vec{\rho})=\int_E x^{\log_q\mu}\D\psi(x) =\int_E
\mu^{\log_q x}\D\psi(x).$$
Letting $\mu\downarrow0$, we obtain
$$I_0(0,\vec{\rho})=\mbox{the jump of $\psi$ at $x=q$,}$$
and this is bounded because $\D\psi$ is a (bounded) Borel measure.
Therefore, the only requirement is $\gamma_k\neq0$, $k\in\Zp$, and we
thus lose no generality by requiring $\gamma_k\equiv1$. 

Since $1+\log_q\mu=\log_q(q\mu)$, we have
\begin{equation}\label{3.1}
I_{k+1}(\mu)=\int_E x^{\log_q \mu}(x-\sigma_{k+1})\rho_k(x)\D\psi(x)
=I_k(q\mu)-\sigma_{k+1} I_k(\mu),\qquad k=0,1,\ldots.
\end{equation}
Recalling that the ratio of $I_{k+1}$ and $I_k$ is a rational $[1/1]$
function, \R{3.1} yields the identity
\begin{equation}\label{3.2}
I_k(q\mu)=\left(\frac{\alpha_k+\beta_k\mu}{1+\delta_k\mu}
+\sigma_{k+1}\right) I_k(\mu),\qquad k=0,1,\ldots.
\end{equation}

We distinguish between the following cases:

\vspace{8pt}
\noindent {\bf Case I} $\alpha_k+\sigma_{k+1}\neq0$.

Let
$$A_k=\alpha_k+\sigma_{k+1}\neq0,\qquad
B_k=-\frac{\beta_k+\sigma_{k+1}\delta_k}{\alpha_k+\sigma_{k+1}}.$$
Hence, \R{3.2} and induction yield
\begin{equation}\label{3.3}
I_k(\mu)=A_k^{-\ell} \frac{(-\delta_k\mu;q)_\ell}{(B_k\mu;q)_\ell}
I_k(q^\ell \mu),\qquad \ell=0,1,\ldots.
\end{equation}
We now let $\ell\rightarrow\infty$. Since $q\in(0,1)$,
$\lim_{\ell\rightarrow\infty} I_k(q^\ell\mu)=I_k(0)$. Therefore
boundedness of $I_k(\mu)$ requires $A_k\equiv1$ and we obtain
\begin{equation}\label{3.4}
I_k(\mu,\vec{\rho})=\frac{(-\delta_k\mu;q)_\infty}{(-(\beta_k+(1-\alpha_k)
\delta_k)\mu;q)_\infty} I_k(0,\vec{\rho}),\quad \sigma_{k+1}=1-\alpha_k\qquad
k=1,2,\ldots.
\end{equation}

We substitute \R{3.4} into \R{3.1}, and this gives
\begin{eqnarray*}
&&\frac{(\delta_{k+1}\mu;q)_\infty}{(-(\beta_{k+1}+(1-\alpha_{k+1})
\delta_{k+1})\mu;q)_\infty}I_{k+1}(0)\\
&=&\left\{\frac{(-\delta_k
q\mu;q)_\infty} {(-(\beta_k+(1-\alpha_k)\delta_k)q\mu;q)_\infty}
-(1-\alpha_k) \frac{(-\delta_k\mu;q)_\infty}
{(-(\beta_k+(1-\alpha_k)\delta_k)\mu; q)_\infty} \right\}I_k(0).
\end{eqnarray*}
Moreover, $I_{k+1}(0)/I_k(0)=\alpha_k$. Thus, $\alpha_k\neq0$,
otherwise $I_{k+1}(0)=0$ $\Rightarrow$ $I_{k+1}\equiv0$, in defiance
of regularity. We deduce thus that
\begin{equation}\label{3.5}
\frac{(-\delta_{k+1}\mu;q)_\infty}{(-(\beta_{k+1}+(1-\alpha_{k+1})
\delta_{k+1})\mu;q)_\infty}=\left(1+\frac{\beta_k\mu}{\alpha_k}\right)
\frac{(-\delta_k
q\mu;q)_\infty}{(-(\beta_k+(1-\alpha_k)\delta_k)\mu;q)_\infty}.
\end{equation}

We distinguish between the following sub-cases, by comparing zeros and
poles on both sides of \R{3.5}

\vspace{6pt}
\noindent{\bf Case I.1} $\alpha_k=1$.

The expression \R{3.5} gives
$$\frac{(-\delta_{k+1}\mu;q)_\infty}{(-(\beta_{k+1}+(1-\alpha_{k+1})
\delta_{k+1}) \mu;q)_\infty}=\frac{(-\delta_k q\mu;q)_\infty}
{(-\beta_k q\mu;q)_\infty}.$$
Therefore $\beta_{k+1}=q(\beta_k-(1-\alpha_{k+1})\delta_k$,
$\delta_{k+1}=q\delta_k$ and $\sigma_{k+1}=0$.

\vspace{6pt}
\noindent{\bf Case I.2} $\beta_k=\alpha_k\delta_k$.

This is inconsistent with regularity, since then
$I_k(\mu)\equiv\alpha_k$.

\vspace{6pt}
\noindent {\bf Case I.3} $\beta_k=\alpha_k\delta_{k+1}$.

We now have
$$\frac{(-\delta_{k+1}q\mu;q)_\infty}{(-(\beta_{k+1}+(1-\alpha_{k+1})
\delta_{k+1})\mu;q)_\infty}=\frac{(-\delta_k q\mu;q)_\infty}
{(-(\beta_k+(1-\alpha_k)\delta_k)\mu;q)_\infty}.$$
Therefore $\delta_{k+1}=\delta_k$ and we are back to the case
I.2, which is irregular.

\vspace{6pt}
Since the above three subcases exhaust all possibilities, we deduce
that only I.1 may happen, hence it {\em must\/} occur for
all $k\in\Zp$. In other words, we have
\begin{equation}\label{3.6}
\left.\begin{array}{rcl}
\sigma_{k+1}&\equiv&0,\\
\alpha_k&\equiv&1,\\
\beta_k&=&q^k\beta_0,\\
\gamma_k&\equiv&1,\\
\delta_k&=&q^k\delta_0,\qquad k=0,1,\ldots.\qquad\qquad
\end{array}\right\}
\end{equation}
Moreover, the moments are
\begin{equation}\label{3.7}
I_k(\mu)=\frac{(-\beta_0\mu;q)_k}{(-\delta_0\mu;q)_k} I_k(0), \qquad
k=0,1,\ldots.
\end{equation}

Let us suppose that $\beta_0\neq0$ and that $\D\psi$ is an atomic
measure with jumps of 
$$\frac{(-\beta_0\mu;q)_\infty}{(-\delta_0\mu;q)_\infty}
\times\frac{\left(\frac{\delta_0}{\beta_0};q\right)_\ell} {(q;q)_\ell}
(-\beta_0)^\ell$$ 
at $q^\ell$, $\ell\in\Zp$. Thus, by the Gau\ss--Heine theorem \cite{S1},
$$I_k(\mu)=\frac{(-\beta_0\mu;q)_\infty}{(-\delta_0\mu;q)_\infty}
{}_1\phi_0 \left[\begin{array}{l}\frac{\delta_0}{\beta_0};
\\\mbox{---}; \end{array} q,-\beta_0q^k\mu\right]=
\frac{(-\beta_0\mu;q)_k} {(-\delta_0\mu;q)_k},\qquad k=0,1,\ldots,$$
where ${}_1\phi_0$ is a generalized basic hypergeometric function.
Hence, the moments are as required by \R{3.7}.

In order for $\D\psi$ to be a Borel measure we require that the jumps
are nonnegative. If $\beta_0>0$ then this implies $\delta_0
q^\ell\geq\beta_0$ for all $\ell\in\Zp$. Therefore $\delta_0>0$ and
only a finite number of jumps might be nonzero. The latter is possible
only if $\delta_0=q^{-N}\beta_0$ for some $N\in\Zp$ and we recover the
{\em $q$-Krawtchouk transformation\/} \cite{IN2}. On the other hand, $\beta_0<0$
and nonnegative jumps can coexist with any $\delta_0>\beta_0$ and we
obtain the {\em Wall transformation\/} \cite{IN2}. The support of the $\mu$s is
the interval $D=(0,|\beta_0|^{-1})$, whereas $E=(0,1)$, thus the
transformation reads
$$T\left\{\sum_{k=0}^m d_k (-\beta_0 x;q)_k (-\delta_0 q^k x;q)_{m-k}
\right\}=\sum_{k=0}^m d_kx^k$$
and it maps polynomials with zeros in $D$ into polynomials with zeros
in $E$ -- compare with $(4.1)$ in \cite{IN2}.

Finally, we consider the case $\beta_0=0$. We derive this as a
limiting case of the Wall distribution. Letting $\beta_0\uparrow0$, we
obtain jumps of
$$\frac{1}{(-\delta_0\mu;q)_\infty} \delta_0^\ell q^{\frac12
(\ell-1)\ell}$$
at $q^\ell$, $\ell\in\Zp$. Thus, $\delta_0>0$ is necessary and
sufficient for $\D\psi$ to be a measure. The transformation becomes
\begin{equation}\label{3.8}
T\left\{\sum_{k=0}^m d_k (-\delta_0 q^k x;q)_{m-k}
\right\}=\sum_{k=0}^m d_kx^k
\end{equation}
and it maps zeros from $[0,\infty)$ to $(0,1)$. We call \R{3.8} -- the
only transformation from the three aforementioned choices of $\omega$
that has been missed in \cite{IN2} -- the {\em Wall$\,{}_0$ transformation.\/}

\vspace{8pt}
\noindent {\bf Theorem 4} The only transformations consistent with the
asserted form of $\D\varphi$ are the Wall, Wall${}_0$ and
$q$-Krawtchouk transformations.






\end{document}